# Barry Mazur contre Pappus

## 15 siècles d'erreur à propos d'une proposition d'Euclide

Salomon Ofman

**Résumé**

L'un des rares points des mathématiques pré-euclidiennes considéré comme établi, est que les bases de la théorie de l'irrationalité développée au livre X des *Éléments* d'Euclide ont été posées par Théétète, un mathématicien grec du 4$^{ème}$ siècle **BCE**. La quasi-totalité des informations sur lui se trouvent dans le dialogue éponyme de Platon, en particulier dans ce que l'on appelle la 'partie mathématique du *Théétète*'. Remontant à Pappus, un mathématicien grec du 4$^{ème}$ siècle **CE**, on considère unanimement que le jeune Théétète y expose une propriété à l'origine de, et généralisée par, la proposition 9 du livre X des *Éléments*.

Il y a quelques années, Barry Mazur a publié un article où il montrait, entre autre, que les deux résultats se rapportaient à des objets mathématiques très différents (« How did Theaetetus proves his Theorem? », p. 234-240). Dans cet exposé, nous allons reprendre cette question en revenant sur le texte du *Théétète*, exposer l'analyse de Mazur et tenter de comprendre comment et pourquoi une telle erreur, '*a strange delusion*' dit Mazur, a pu s'imposer aux historiens et aux philosophes, sur une aussi longue période.

Cet exposé est basé sur un travail en commun avec Luc Brisson (CNRS-Centre Jean Pépin, ENS) sur l'interprétation de la 'partie mathématique du *Théétète*'.

## I. Présentation du 'passage mathématique du *Théétète*'

Très peu de choses sont considérées établies de manière fiable concernant les mathématiques avant Euclide (vers le 3$^{ème}$ siècle BCE) par les historiens des mathématiques contemporains.

[Dhombres - Avant_Euclide.mp4](Dhombres - Avant_Euclide.mp4)

Il y a deux ans, j'avais déjà commis ici un exposé sur ce passage fondamental aussi bien pour l'histoire des mathématiques, que pour le dialogue lui-même, et plus généralement la philosophie de Platon où, comme on le sait, les mathématiques jouent un rôle crucial. Pour les historiens des mathématiques, son importance tient en ce qu'il est le seul texte de cette époque donnant des informations sur les mathématiques pré-euclidiennes, en particulier sur l'origine de la théorie de l'irrationalité développée dans les *Éléments*. Je vais brièvement en rappeler les points essentiels.

- L'objet du *Théétète* est de définir ce qu'est la science ('*epistêmê*'). Le principal texte concernant les mathématiques pré-euclidiennes se trouve dans ce que l'on appelle 'la partie mathématique du *Théétète*' (147d-148b).
- Dans une première partie, Théétète un très jeune Athénien présenté comme extrêmement brillant, rend compte d'une leçon sur les 'puissances' donnée par un mathématicien déjà âgé et reconnu, Théodore de Cyrène. Le sujet en est en termes modernes l'étude de l'irrationalité des racines carrée des entiers 3, 5, jusqu'à 17. Il y a un nombre considérable de travaux sur l'interprétation de cette leçon et aussi bien que sur le sens du texte dont quasiment aucun mot n'est pas objet de désaccord parmi les chercheurs modernes. Un point important à retenir est que la démonstration de Théodore est faite au cas par cas, et ne résulte pas d'un résultat général sur les racines carrées des entiers. Suite à cette leçon, Théétète et un de ses camarades, ayant le même nom que le philosophe Socrate, retravaillent sur cette question pour mieux la comprendre, comme le ferait d'ailleurs aujourd'hui encore de bons élèves. Ils décident de s'atteler à une généralisation pour éviter la procédure au cas par cas, qui est incompatible avec la nature des entiers, en tant que multitude infinie.

- Nous avons présenté pour être plus facilement compréhensible le problème dans un langage moderne. Toutefois, pour comprendre l'enjeu, il faut se replacer dans le cadre ses mathématiques grecques de l'Antiquité. Trois points sont à garder à l'esprit :

    i) Contrairement à la conception moderne, les entiers ne sont pas considérés comme une partie d'un ensemble plus grand, pour nous les nombres réels. En fait on a deux parties des mathématiques. L'une s'occupe des entiers exclusivement, c'est ce qu'on appelle au temps de Socrate 'la théorie du pair et de l'impair', l'autre des grandeurs géométriques. L'une est encore caractérisée comme travaillant sur le discret, l'autre le continu. Si ces deux parties des mathématiques sont bien différenciées, il n'empêche qu'on peut représenter géométrique les entiers par des grandeurs géométriques, en posant un segment de droite arbitraire comme représentant l'unité arithmétique. C'est pourquoi, la géométrie va occuper une place privilégiée dans les mathématiques grecques.

    ii) Dans la représentation moderne, on considère les nombres formant une suite croissante partant des entiers positifs (avec toutefois 0), puis les rationnels et enfin les réels (pour simplifier, on ne considère pas les nombres négatifs).
    En mathématique grecque, on a tout d'abord les entiers positifs (sans 0, l'unité étant considérée suivant les cas comme un nombre ou pas). Puis on trouve les grandeurs (géométriques) irrationnelles, ou plutôt incommensurables, le terme irrationnels très souvent utilisé dans le *Théétète*, n'apparaît pas dans le cadre mathématique. L'incommensurabilité est une propriété relative supposant non pas une mais deux grandeurs. Comme son nom l'indique, deux grandeurs géométriques *a* et *b* sont incommensurables s'il n'existe pas de mesure commune aux deux. Autrement dit, il n'existe pas de grandeurs *u* telle *a* et *b* soit simultanément des multiples entiers de *u*.
    Il est facile de réécrire cela en langage moderne : *a* et *b* sont commensurables si l'on peut écrire $a/b = mu/nu = m/n$, et incommensurables dans le cas contraire. En d'autres termes, *a* et *b* son commensurables si et seulement si $a/b$ est un nombre rationnel, incommensurable si $a/b$ est irrationnel.
    La connexion entre commensurabilité/incommensurabilité et rationalité/irrationalité se fait de en posant une grandeur arbitraire comme étant l'unité. Une grandeur est alors rationnelle ou irrationnelle si elle est commensurable ou incommensurable à cette unité.

iii) Ceci paraît fort simple. Il y a toutefois un problème essentiel, c'est que la notion de *nombre* **rationnel** n'apparaît pas dans les mathématiques grecques. Il n'est donc pas possible de parler du *nombre a/b*. De fait, on a là un nouvel objet mathématique qui est le *rapport* des nombres *a* et *b*. Mais le *rapport* n'est pas un nombre, même s'il est souvent confondu avec la fraction *a/b*. Jusqu'ici on n'est pas si loin du traitement moderne. Les choses sont différentes lorsqu'on passe aux grandeurs géométriques, en particulier lorsque ces grandeurs sont incommensurables. Le problème est qu'on n'a pas l'ensemble des nombres réels qui permet de travailler indépendamment de savoir si l'on a affaire à des grandeurs commensurables ou incommensurables. C'est le livre V des Éléments qui permet de pallier à cette absence, et d'obtenir une théorie générale des rapports, que l'on ait des entiers ou des irrationnels. Mais encore une fois, il n'est pas question de *nombres* rationnels, mais de *rapports* commensurables ou incommensurables.

a) La méthode de Théodore de construction des grandeurs irrationnelles. Bien que cela ne soit pas explicite, il est raisonnable de penser qu'elle se fonde sur le corollaire immédiat suivant du théorème de Pythagore, que dans un triangle rectangle OBD, on a :

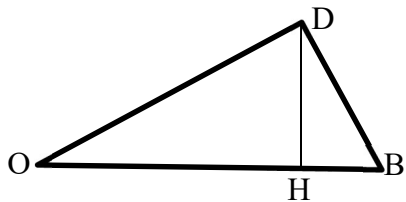

*Figure 1*

Le carré de côté OD est égal (en surface) au rectangle de côtés OH et HB.
D'où la construction des grandeurs irrationnelles que nous appelées, les racines carrées de 3, 5, …, 17 :

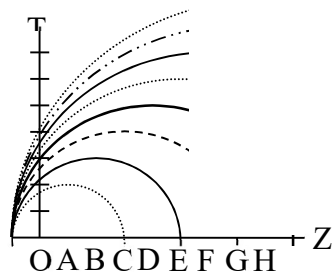

*Figure 2*

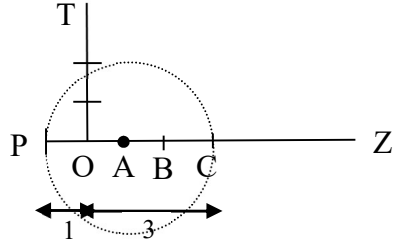

*Figure 3-Détail : Le premier cas, l'entier 3*

## II. Le texte du *Théétète*

Pour les curieux, nous donnons le texte original de la partie qui nous concerne ici, et la traduction que nous en avons faite avec Luc Brisson.

| THEAETETUS<br>We divided the integer in its totality in two. One which has the power to be the product of an equal integer by itself [**148a**] we likened to the square as a figure, and we called it a square or an equilateral integer.<br><br>SOCRATES<br>Good, so far.<br><br>THEAETETUS<br>Any integer found between two [successive square or equilateral integers], such as three or five or any integer that can't come from the product of an equal integer by an equal integer, but that comes from the product of unequal integers, a greater and a smaller, is always encompassed by a longer and a less side, having depicted it in shape as rectangular figure, we called it a rectangular integer.<br><br>SOCRATES<br>That's excellent. But how did you go on?<br>THEAETETUS<br>All the segments of line which form the sides of the figure squaring an equilateral and plane integer, we defined as "length", while we defined as "powers" all the segments of line which form the sides of the figure squaring a rectangular integer, because, although they are not commensurable as length with the formers, [**148b**] they are as | ΘΕΑΙ.<br>Τὸν ἀριθμὸν πάντα δίχα διελάβομεν· τὸν μὲν δυνάμενον ἴσον ἰσάκις γίγνεσθαι τῷ τετραγώνῳ τὸ σχῆμα ἀπεικάσαντες τετράγωνόν τε καὶ ἰσόπλευρον προσείπομεν.<br><br>ΣΩ.<br>Καὶ εὖ γε.<br><br>ΘΕΑΙ.<br>Τὸν τοίνυν μεταξὺ τούτου, ὧν καὶ τὰ τρία καὶ τὰ πέντε καὶ πᾶς ὃς ἀδύνατος ἴσος ἰσάκις γενέσθαι, ἀλλ' ἢ πλείων ἐλαττονάκις ἢ ἐλάττων πλεονάκις γίγνεται, μείζων δὲ καὶ ἐλάττων ἀεὶ πλευρὰ αὐτὸν περιλαμβάνει, τῷ προμήκει αὖ σχήματι ἀπεικάσαντες προμήκη ἀριθμὸν ἐκαλέσαμεν.<br><br>ΣΩ.<br>Κάλλιστα. ἀλλὰ τί τὸ μετὰ τοῦτο;<br>ΘΕΑΙ.<br>Ὅσαι μὲν γραμμαὶ τὸν ἰσόπλευρον καὶ ἐπίπεδον ἀριθμὸν τετραγωνίζουσι, μῆκος ὡρισάμεθα, ὅσαι δὲ τὸν ἑτερομήκη, δυνάμεις, ὡς μήκει μὲν οὐ συμμέτρους ἐκείναις, τοῖς δ'ἐπιπέδοις ἃ δύνανται. καὶ περὶ τὰ στερεὰ ἄλλο τοιοῦτον. |

| areas they have the power to produce. And it is the same in the case of solids. | |

Ce passage pose des difficultés de traduction et d'interprétation redoutables. Toutefois son contenu purement mathématique est relativement aisé à comprendre.

### i) L'argument de Théétète

La construction et la difficulté soulevée est transparente pour les modernes. En effet, le problème difficile ici et de concilier nombre et figures, arithmétique et géométrie. Théétète en effet dans un premier temps explique comment par une représentation figurée, on peut unifier géométriquement entiers et grandeurs y compris grandeurs incommensurables.

**1) Représentation géométrique des nombres**

Tout d'abord, il considère une représentation des entiers, non pas par une ligne comme chez Euclide, et certainement déjà au temps de Socrate, mais par un rectangle de côtés entiers. Tout entier, dit-il est en effet représentable par un certain rectangle dont la surface est égale à cet entier. Deux choses ne sont pas explicitées.
   a) Nécessairement l'unité doit être alors considérée comme un nombre entier. En effet, la représentation de tout nombre premier se fait uniquement par un rectangle dont l'un des côtés est l'unité et l'autre ce nombre lui-même.
   b) Cette représentation n'est pas unique. En effet, si tout entier peut être représenté à la manière des nombres premiers, pour les autres il existe diverses possibilités. Ainsi, et c'est le cas fondamental, si un entier est un carré parfait, il est représentable non pas comme un simple rectangle, mais bien comme un carré.

**2) Les quadratures entières**

Ensuite, il découpe l'ensemble des nombres entre ceux qui sont précisément représentables par des carrés (les carrés parfaits) et tous les autres. Alors,
   a) Pour ce qui est des carrés parfaits, les côtés du carré seront entiers.
   b) Par contre pour ce qui est des autres entiers (non carrés parfaits), on peut associer au rectangle les représentants un carré, mais cette fois ce côté non entier.
   c) À chaque rectangle précédent, on peut associer un carré de même surface, par la méthode de la figure 2 ci-dessus.

On a ainsi associé à chaque les entier, un certain carré qui en langage moderne a pour côté la racine carrée de cet entier (mais cela n'est évidemment pas présenté ainsi ici).

**3) Généralisation de la leçon de Théodore**

Ce que désire ici Théétète, c'est généraliser à tous les entiers, les résultats de Théodore arrêté à 17. Cela est dit clairement en introduction à cette leçon qui illustre la signification d'une définition, comme universel à la totalité d'une classe généralement infinie (ici les entiers). En termes modernes, Théodore a montré que la racine carrée d'un entier était ou bien un entier, ou bien irrationnel.
De même en termes modernes, Théétète a procédé ainsi :
Il a divisé l'ensemble $N$ des entiers (privé de 0) en deux sous-ensembles $P$ et $R$, respectivement les entiers qui sont des carrés et ceux qui ne le sont pas. Puis il a considéré les ensembles ***P'*** et ***R'*** formés respectivement des racines carrées des entiers de ***P*** et de ***R***.

D'où dit-il, on obtient
- D'une part, des entiers, donc *évidemment commensurables* à l'unité et entre eux, ce sont les éléments de ***P'***.
- D'autre part, des éléments *tout aussi évidemment incommensurables* à l'unité et aux éléments de ***P'***.

**Une** (parmi d'autres) traduction en langage moderne du résultat énoncé par Théétète est la suivante :

**Proposition A**. *La racine carrée d'un entier est rationnelle si et seulement si cet entier est un carré parfait.*

Comparé au récit de Théétète, on voit immédiatement de la concision moderne, là où un jeune grec a besoin d'une demi-page.

### III. Le texte de Mazur

> "*As I have already mentioned, we will be discussing the proof of this theorem, as it appears in the extant ancient literature. What is strange, though, is that a popular delusion seems to be lurking in the secondary literature on this topic. Specifically, you will find – in various places – the claim that Theaetetus' theorem is proven in Proposition 9 of Book X of Euclid's Elements. It doesn't serve any purpose here to list the places where you find this incorrect assertion, except to say that it is incorrect, and it remains a thriving delusion since at least one important article published as late as 2005 repeats it. It is an especially strange delusion since nothing subtle is going on here. Even a cursory glance at Proposition 9 of Book X will convince you that what is being demonstrated there – if you take it in a modern perspective - is an utter triviality. Proposition 9 of Book X stands, though, for an important issue in ancient thought if taken on its own terms, but it won't prove irrationality of anything for us, let alone irrationality of all the numbers that Theaetetus proves. One might imagine that Heath's commentary on this – which is perfectly clear, and says exactly what is indeed proved in Proposition9 – would dispel the misconception that Theaetetus' theorem about the irrationality of surds is contained in this proposition, but it seems that this has held on with some tenacity. I would guess that the source of this error is quite early, as early as the commentaries of Pappus, but I offer this guess timidly because that would seem to imply that poor Proposition9 of Book X has been often cited but far less often read with attention since the fourth century AD.*" (p. 235).

### IV. La proposition X.9 des *Éléments* d'Euclide

Voici la traduction de Thomas Heath :

**Proposition B** (proposition X.9 des *Éléments*). *The squares on straight lines commensurable in length have to one another the ratio which a square number has to a square number; and squares which have to one another the ratio which a square number has to a square number will also have their sides commensurable in length. But the squares on straight lines incommensurable in length have not to one another the ratio which a square number has to a square number, and squares which have not to one another the ratio which a square number has to a square number will not have their sides commensurable in length either.*

Voici **une** (parmi d'autres) traduction en langage moderne de cette proposition :

**Proposition B**. *Le carré d'un nombre est le carré d'un rationnel si et seulement si ce nombre est rationnel. Et inversement, le carré d'un nombre n'est pas le carré d'un rationnel si et seulement si ce nombre n'est pas rationnel.*

On notera là encore, la différence de longueur entre les deux énoncés, qui n'est certainement pas dû au langage grec en soi, qui est plutôt concis (voir le texte ci-dessus du *Théétète* et comparez-le à sa traduction).

**V.    *Théétète* et La proposition X.9 d'après Pappus**

Une remarque préliminaire. Le texte ne dit pas de quelle proposition des Éléments il s'agit, mais d'après ce les propriétés discutées, c'est très probablement de la proposition 9 qu'il est question.

**1) Voici comment Pappus présente le récit donné dans le *Théétète***

> Theodorus had discussed with Theaetetus the proofs of the *powers* (i. e. squares)[66] which are commensurable and incommensurable in length relatively to the *power* (square) whose measure is a [square] foot[67], the latter had recourse to a general definition of these *powers* (squares), after the fashion of one who has applied himself to that knowledge which is in its nature certain (or exact)[68]. Accordingly he divided all numbers into two classes[69]; such as are the product of equal sides (i. e. factors)[70], on the one hand, and on the other, such as are contained by a greater side (factor) and a less; and he represented the first [class] by a square figure and the second by an oblong, and Ms. 26 r concluded that the *powers* (squares) which *square* (i. e. form into a square figure) a number whose sides (factors) are equal[71], are commensurable both in square and in length, but that those which *square* (i. e. form into a square figure) an oblong number, are incommensurable with the first [class] in the latter respect (i. e. in length), but are commensurable occasionally with one another in one respect[72].

Il semble que Pappus (ou l'un des divers rédacteurs et traducteurs) ne lit pas directement le texte de Platon, mais à travers un commentaire. En effet,

- d'une part, on dit qu'il s'agit d'une discussion entre Théodore et Théétète, alors que le texte parle d'une leçon donnée par Théodore.
- D'autre part, on dit que les 'puissances' correspondant aux nombres qui ne sont pas des carrés parfaits (i.e. leurs racines carrées) sont incommensurables aux puissances correspondant aux nombres carrés parfaits (i.e. les entiers), ce qui est bien dans le récit, mais également qu'ils sont parfois commensurables entre eux, ce qui n'est indiqué d'aucune façon dans le texte de Platon.
- Ici, on présente Théétète seul qui suivant les règles de la science la plus exacte, donne une définition générale, là où son professeur, Théodore, s'est borné à parler de grandeurs particulières. Dans le dialogue platonicien, Théétète insiste qu'il n'a pas travaillé seul, mais en partenariat avec un de ses camarades, nommé comme Socrate, comme le philosophe à qui il s'adresse.

Un point mérite attention dans le passage ci-dessus. Il apparaît que d'une part Théodore, d'autre part Euclide, *prouvent* leurs résultats. Quant à Théétète, il est dit simplement qu'il *conclut* de sa construction. Ne connaissant pas l'arabe, je ne sais pas s'il faut en tirer quelque chose.

## 2) La proposition X.9 selon Pappus

Voici la suite du texte précédent qui concerne cette fois Euclide :

> another in one respect[72]. Euclid, on the other hand, after he had examined this treatise (or theorem) carefully for some time and had determined the lines which are commensurable in length and square, those, namely, whose *powers* (squares) have to one-another the ratio of a square number to a square number, proved that all lines of this kind are always commensurable in length[73]. The difference between Euclid's statement (or proposition)[74] and that of Theaetetus which precedes it, has not escaped us. The idea of determining these *powers* (squares) by means of the square numbers is a different idea altogether from that of their having to one-another the ratio of a square [number] to a square [number][75]. For example, if there be taken, on the one hand, a *power* (square) whose measure is eighteen [square] feet, and on the other hand, another *power* (square) whose measure is eight [square] feet, it is quite clear that the one [power or square] has to the other the ratio of a square number to a square number,

Page 11.

— 74 —

> the numbers, namely, which these two double[76], notwithstanding the fact that the two [powers or squares] are determined by means of oblong numbers. Their sides, therefore, are commensurable according to the definition (thesis) of Euclid, whereas according to the definition (thesis) of Theaetetus they are excluded from this category. For the two [powers or squares] do not *square* (i. e. do not form into a square figure) a number whose sides (factors) are equal, but only an oblong number. So

L'opposition entre Théétète et Euclide, outre l'insistance sur le caractère de preuve mise en œuvre par ce dernier, est la différence de leurs constructions. En effet, alors que Théétète ne s'intéresse qu'aux entiers ('déterminer les puissances par les entiers carrés') alors qu'Euclide travaille directement sur des rapports commensurables ('ayant l'un l'autre un rapport d'un entier carré à un entier carré'), ce qui correspond en termes modernes aux rationnels.

Cette manière de procéder est très différente et n'a pas échappé à Pappus, dit-il, qui insiste sur la bien plus grande généralité de la définition d'Euclide. Il en donne un exemple. Alors que dans le cadre de Théétète, concernant exclusivement les entiers, les grandeurs telles que les racines carrées de *8* et de *18* échappent complètement, pour Euclide elles sont commensurables, puisque *√18/√8 = √9/√4 = 3/2*.

Autrement dit, le résultat d'Euclide est bien plus général que celui de Théétète.

## VI. Le sens de la proposition X.9

Voici une 'démonstration' moderne la proposition B.

**Proposition B**. *Le carré d'un nombre est le carré d'un rationnel si et seulement si ce nombre est rationnel. Et inversement, le carré d'un nombre n'est pas le carré d'un rationnel si et seulement si ce nombre n'est pas rationnel.*

**Démonstration.** Soient *t* une grandeur et *m, n* des entiers. On a alors :

$t^2 = (m/n)^2 = m^2/n^2 \Leftrightarrow t = m/n$. Quant à la seconde partie, ce n'est que la contraposée de la première.

La démonstration dans les *Éléments* ne nécessite pas moins d'une page et demi. Pourquoi une telle différence de traitement entre la démonstration des *Éléments* et celle algébrique moderne ?

- D'une part, dans la seconde on traite simultanément des nombres rationnels et irrationnels (sous-entendu, en tant qu'élément des réels), alors que dans celle d'Euclide, on doit considérer et comparer des nombres et des grandeurs géométriques. On ne peut pas écrire *t*, mais considérer le rapport de deux segments de longueurs respectives *a* et *b*.
- L'autre difficulté est que dans les *Éléments*, on n'a pas affaire à des nombres, mais à des rapports. Ainsi la dernière égalité dans la démonstration moderne $(m/n)^2 = m^2/n^2$ est loin d'être évidente, car il n'est déjà pas clair ce que signifie le carré d'un rapport, moins encore qu'il est 'égal' au rapport des carrés. Dans le cadre des mathématiques grecques, **la proposition X.9 est donc très loin d'être triviale,** et bien au contraire démontre une propriété fondamentale pour la théorie de l'irrationalité.

## VII. La proposition A du *Théétète*.

Nous allons maintenant considérer la proposition du Théétète, elle aussi sous sa forme modernisée.

**Proposition A**. *La racine carrée d'un entier est rationnelle si et seulement si cet entier est un carré parfait.*

**Proposition A**. *Un entier est un carré parfait si et seulement si sa racine carrée est rationnelle.*

Qu'en est-il de sa démonstration algébrique moderne. Il s'agit de montrer que pour tout entier *r* :

$r = p^2 \Leftrightarrow \sqrt{r} = m/n$ (*m* et *n* entiers).

Le début de la démonstration est la même que celle de la proposition B :

$r = (m/p)^2 = m^2/p^2 \Leftrightarrow \sqrt{r} = m/n$,

où *r* joue le rôle de $t^2$.

Une démonstration moderne passe par le *théorème fondamental de l'arithmétique*, i.e. la décomposition unique (à permutations près) de tout entier comme puissances de nombres premiers. La dernière égalité donne en effet : $rn^2 = m^2$, d'où l'unicité de la décomposition implique que $r$ soit un produit de puissances paires d'entiers, donc de carrés, donc un carré lui-même.

On va la donner dans un cadre euclidien, mais sous une écriture modernisée. Tout d'abord dans l'écriture :

$\sqrt{r} = m/n \Leftrightarrow r = (m/n)^2 = m^2/n^2$,

on peut choisir $m$ et $n$ **premiers** entre eux. Cela est l'objet d'une proposition d'Euclide, où il montre

- que ce sont les entiers minimaux donnant cette fraction (proposition 22) d'où la représentation d'un rapport par des éléments premiers est unique,
- et surtout que si deux entiers $a$ et $b$ ont même rapport, alors $m$ divise $a$ et $n$ divise $b$ (proposition 20). Cela résulte de la définition 20 définissant le rapport de deux entiers, et de la possibilité de permuter les éléments d'un rapport (proposition 13).

Les entiers $m$ et $n$ étant choisis relativement premiers, leurs carrés $m^2$ et $n^2$ le sont également (proposition VII. 24). L'égalité $r = m^2/n^2$ s'écrit sous forme de rapport : $r/1 = m^2/n^2$ et par unicité (due à la minimalité), on a : $r = m^2$ (et $n^2 = 1$) et $r$ est un carré parfait.

Le passage de de deux nombres premiers entre eux à leurs carrés est essentiellement ce qu'on appelle aujourd'hui le lemme de Gauss.

**VIII. Le texte du *Théétète* et la proposition X.9 : un double paradoxe.**

    **1) Le paradoxe du *Théétète***
    Si l'on retourne au texte du *Théétète*, la seule démonstration est un découpage des entiers permettant leur représentation géométrique. Les seules propriétés utilisées sont triviales, au moins pour nous, et aucun appel à une quelconque forme du lemme de Gauss.

    **2) Le paradoxe de la proposition X.9**
    Cette fois, en comparant les démonstrations modernes de la proposition X.9 et de la proposition A, l'une résulte d'une équivalence immédiate, l'autre nécessite un résultat profond d'arithmétique, qui d'ailleurs est spécifique aux entiers (et plus généralement lorsqu'on se place dans un anneau factoriel i.e. lorsque tout élément admet une décomposition unique, à permutation près, en facteurs premiers).

Sauf à penser que le passage à l'arabe ou au grec ancien transforme la difficulté d'un résultat mathématique, il semble donc qu'on ait un miracle, à savoir une manière triviale de démontrer un résultat qui ne l'est pas (de notre point de vue toujours).

Cela est évidemment choquant pour un mathématicien qui ne croit pas aux miracles mathématiques du moins. C'est aussi le second paradoxe qui avait choqué Mazur, qui ne s'intéressait pas particulièrement à la question du *Théétète*, mais bien plutôt à la comparaison entre la proposition X.9 et la proposition A.

## IX. La propositions X.9 et de la proposition A

Nous allons tout d'abord considérer le second paradoxe.

Tout d'abord, si en effet dans les *Éléments*, les deux propositions apparaissent bien, c'est dans des livres différents. L'une dans le cadre de géométrique des grandeurs irrationnelles du livre X, l'autre dans le cadre arithmétique du livre VII.

Reprenons donc la question en termes modernes.

- La proposition X.9 affirme que la racine carrée d'un nombre rationnel est rationnelle si et seulement si ce nombre est le carré d'un nombre rationnel.
- La proposition A dit que la racine carrée d'un entier est rationnelle si et seulement si ce nombre est le carré d'un entier.
- Est-il possible que malgré la différence de difficulté dans les démonstrations, l'une soit la généralisation de l'autre ? Cela n'est absolument pas le cas. En effet, soit la racine carrée d'un entier $r$, donc d'un nombre rationnel ; d'après l'équivalence de la proposition X.9, si cette racine est rationnelle, alors cet entier est le carré d'un *nombre rationnel*. Mais rien ne dit, et c'est précisément là que se cache la difficulté, que ce soit le carré d'un nombre entier, ce qui est pourtant la conclusion de la proposition A. Autrement dit, $r$ est sûrement le carré d'un nombre rationnel, mais rien ne prouve qu'il soit un carré parfait !

## X. L'énoncé du *Théétète* et la proposition A

Considérons alors le premier paradoxe, la conclusion du *Théétète* donnant la proposition A, à partir d'un découpage des entiers. Là encore, pour les modernes, on est dans le cadre d'une pure trivialité face à une propriété passant par un résultat fondamental, le théorème fondamental de l'algèbre.

Le problème ici étant non pas que les énoncés sont différents, la proposition A étant précisément cet énoncé, mais dans la preuve (ou absence de preuve). Le texte de Platon est particulièrement 'pervers' dans le sens où la propriété essentielle est cachée, et on passe sur elle comme une évidence. Mais essayons d'abord de comprendre ce que Théétète peut bien vouloir dire, suite à la leçon mathématique qu'il a suivie. Là encore, par commodité, on utilisera le langage moderne, tout en gardant en tête les dangers que cela comporte.

Dans cette leçon, Théodore considère donc des racines carrées d'entiers, 3, 5,... et jusqu'à 17. Quels sont exactement ces nombres est très discuté, pour faire simple, ce sont les impairs compris de 3 à 17, 3 inclus et 17 exclu. Comme on le voit dans cet échantillon deux cas sont possibles : ou bien la racine carrée est irrationnelle, ou bien elle est entière (cas de 9). Par ailleurs, l'objet de la leçon est de montrer l'irrationalité de certaines grandeurs (les 'puissances'), quant à Théétète et son camarade, il veut donner une définition générale qui comprend la totalité des 'puissances' autrement dit toutes les racines carrées irrationnelles. C'est l'objet de sa division des entiers entre 'carrés parfaits' et 'non carrés', la racine carrée des premiers donnant des entiers, celle des seconds, les 'puissances', c'est-à-dire les irrationnels. C'est bien d'ailleurs ce qu'il affirme, car si l'on reprend sous une forme moderne ce qu'il dit à ce propos :

'Les racines carrées les nombres non carrés seront définis comme les 'puissances' parce que, bien qu'elles ne soient *pas commensurables* avec les entiers, leurs carrés le sont.'

Les 'puissances', c'est-à-dire ce que Théétète et son camarade cherchent à définir, suite à la leçon de Théodore, sont des irrationnelles, mais de carré rationnel. Le caractère d'irrationalité est donc pour eux une évidence et le résultat qu'ils ont donc en tête est le suivant :

**Proposition A' de Théétète** : *la racine carrée d'un entier est un entier si et seulement si cet entier est un carré.*

Ce résultat est en effet une évidence, car il s'écrit : pour tout entier $m$,

$\sqrt{m} = n \Leftrightarrow m = n^2$ (où $n$ est un entier).

Ce qui échappe précisément à Théétète et à son camarade, c'est la difficulté à prouver que la racine carrée d'un entier non carrée peut ne pas être entière sans être irrationnelle. C'est sans doute dans la littérature, le premier exemple humien en mathématique : l'impossibilité de généraliser une propriété même si elle est démontrée dans un grand nombre de cas singuliers.

XI. **Retour à Pappus et à une erreur de 15 siècles**

Dans ce passage, dit Théétète, il s'agit d'une construction dont le but est d'obtenir les 'puissances', c'est-à-dire des grandeurs incommensurables dont le carré est commensurable. C'est bien ce qu'un mathématicien qui lit ce passage va retenir. D'après le texte, il y a implicitement que ce sont les racines carrées des entiers, ce qu'un mathématicien après Euclide, et même sans doute après Théétète, sait être vrai. Mais il sait aussi que ce ne sont pas les seules, d'où la supériorité d'Euclide soulignée par Pappus, qui élargit considérablement cet ensemble. Au lieu des seules racines carrées d'entiers non carrés, ce sont toutes celles qui ne sont pas des carrés de grandeurs rationnelles. Pappus n'a donc pas de raison impérieuse d'examiner comment sont obtenus ces résultats du point de vue mathématique, les dialogues de Platon ne sont pas des traités de mathématique. Les démonstrations se trouvent dans le livre X des Éléments, si comme le dit la tradition, Théétète est à l'origine, voire l'auteur d'une partie de ce livre.

Si comme nous l'avons dit, Pappus ne lit pas directement Platon, mais un commentaire portant sur le travail de Théétète, l'erreur sur la relation entre le passage étudié et la proposition X.9 doit être très antérieure et reprendre une tradition déjà bien établie au 5[ème] siècle CE. Toutefois, pour les mathématiciens qui ne sont concernés que par l'aspect historique du texte concernant l'origine de la théorie de l'irrationalité du livre X et non son caractère mathématique, erreur ou pas, cela ne leur importe pas vraiment. Avec ou sans démonstration, Théétète donne, comme le dit Pappus, un ensemble de grandeurs 'incommensurables en longueur' mais 'commensurables en carré', qu'Euclide va généraliser. Il en va différemment pour les exégètes de Platon, car il y a un écart considérable entre ce qu'il dit faire, donner une définition, et ce qu'il fait réellement, en énonçant une propriété fondamentale des entiers. Pour les modernes, le problème n'est pas son rapport à la proposition X.9, mais de donner un sens à un énoncé vrai, mais dont la démonstration non seulement manque, mais semble inutile. Erreur qui est néanmoins commune en mathématiques, croire une propriété évidente alors qu'elle est au contraire très difficile. Pour prendre un exemple contemporain, cela est apparu plusieurs fois au cours de la démonstration du 'grand théorème de Fermat', erreurs successivement corrigées pour aboutir à une démonstration satisfaisante pour les mathématiciens. Une fois cela établi, il reste à comprendre pourquoi Platon présente cette situation dans le cadre d'une recherche d'une définition de la science, mais ce n'est plus du domaine d'étude de Mazur et donc de cet exposé.

## Brief Bibliography